\documentstyle{article}

\newtheorem{theorem}{Theorem}[section]
\newtheorem{lemma}{Lemma}[section]
\newtheorem{corollary}{Corollary}[section]
\newtheorem{definition}{Definition}[section]
\newtheorem{proposition}{Proposition}[section]

\begin{document}

\title{Lyapunov functions: necessary and sufficient conditions for practical
stability}
\author{F.G. Garashchenko, O.M.Bashniakov and V.V.Pichkur \\
Kyiv Shevchenko University, Glushkova Prosp., 6, Kyiv,\\
252127, Ukraine. \\
E-mail: chaton@akcecc.kiev.ua, bigcat@briefcase.com}
\maketitle

\section{Introduction}

In many practical problems, it is important to construct the set of all
possible initial states of a dynamical system for which the trajectories do
not violate given phase restrictions. For example, one of the basic problems
in acceleration technology concerns maximal capture of particles. It is
equivalent to maximization of the domain of stability \cite{kniga1}. On the
other hand, the practical stability of motion generalizes the well-known
Chetayev's definition of $\left\{ \lambda ,A,t_0,T\right\} $-stability \cite
{Chetayev}. Analysis of stability on a finite time interval is realized
using Lyapunov functions method. In works \cite{kniga1},\cite{Martinuk},\cite
{IEEE},\cite{Control} have been proved the sufficient conditions for
practical stability and, for example, in \cite{kniga1},\cite{stat1},\cite
{stat2},\cite{stat3} have been obtained the necessary conditions in
particular cases. The similar results hold for parametric and discrete
systems \cite{kniga1},\cite{Pant}.

The purpose of this paper is twofold. First, we shall prove the necessary
and sufficient conditions for practical stability of nonlinear dynamical
system at general phase restrictions. It will be shown, that in such a case
the Lyapunov function is nondifferentiable. But if the set of initial data
is starry compact, then it is possible building Lyapunov function which
belongs to differentiable functions class. Our second goal is to estimate
the optimal sets of initial conditions in structural forms for linear system
and concrete phase restrictions. This estimations are important in
theoretical sense and enable numerical methods to be devised \cite{stat2}.

\section{Auxiliary propositions}

Let $\rho $ be a metric in the space $R^n$. We shall denote by $\partial K$
the set $K$ frontier and by $intK$ the interior of $K$, $K\subset R^n$. The
following proposition is true.

\begin{lemma}
\label{lem} If $K\subset R^n$ is a compact then there exists a function $%
\alpha :R^n\rightarrow R^1$, $\alpha \in C\left( R^n\right) $, such that: $%
\alpha (x)=1$ at $x\in \partial K$; $\alpha \left( x\right) <1$ at $x\in intK
$; $\alpha \left( x\right) >1$ at $x\in R^n\backslash intK$.
\end{lemma}

{\bf Proof.} As it is known, the function $\psi \left( x\right)
=\min\limits_{y\in \partial K}\rho \left( x,y\right) $ is continuous.
Moreover $\psi \left( x\right) \geq 0$ at $x\in R^n$ and $\psi \left(
x\right) =0$ if and only if $x$ belongs to $\partial K$. Consider the
function $\alpha _1\left( x\right) =1-\psi \left( x\right) $, $x\in K$. This
function is continuous on $K$, $\alpha _1\left( x\right) <1$ if $x\in intK$
and $\alpha _1\left( x\right) =1$ if $x\in \partial K$. On the other hand
the function $\alpha _2\left( x\right) =1+\psi \left( x\right) $ belongs to $%
C\left( R^n\backslash intK\right) $, $\alpha _2\left( x\right) >1$, $x\in
R^n\backslash K$ and $\alpha _2\left( x\right) =1$ if $x\in \partial K$.
Then the function $\alpha \left( x\right) =\left\{ 
\begin{array}{c}
\alpha _1\left( x\right) =1-\psi \left( x\right) ,x\in K, \\ 
\alpha _2\left( x\right) =1+\psi \left( x\right) ,x\in R^n\backslash K
\end{array}
\right. $ satisfies the lemma's conditions.

Denote $\left\| x\right\| $ the Euclidean norm of vector $x\in R^n$, $%
S_r\left( a\right) =\{x\in R^n:\left\| x-a\right\| =r\}$, $K_r\left(
a\right) =\{x\in R^n:\left\| x-a\right\| \leq r\}$, $L_x=\{u\in
R^n:u=\lambda x,\lambda \geq 0\}$, $x\in R^n\backslash \left\{ 0\right\} $.

\begin{definition}
\label{starry}The set $E\subseteq R^n$ is said to be starry if it is closed, 
$0\in intE$ and for arbitrary $x\in R^n\backslash \left\{ 0\right\} $ either 
$L_x\subset E$ or there exists a unique point $z\in L_x$ so that $z\in
\partial E$.
\end{definition}

The starry sets are singly connected. In fact, if $A$ is a starry set, $a\in
A$ and $b\in A$ then it is possible to build the curve $z\left( \lambda
\right) =\left\{ 
\begin{array}{c}
\lambda a,\lambda \in \left[ 0,1\right] , \\ 
-\lambda b,\lambda \in \left[ -1,0\right]
\end{array}
\right. $ and $z\left( \lambda \right) \in A$, $\lambda \in \left[
-1,1\right] $, $z\left( 1\right) =a$, $z\left( -1\right) =b$.

If $E$ is a starry compact then: a) for any $l\in S_1\left( 0\right) $ there
exists a unique $\alpha >0$ such that $\alpha l\in \partial E$ and $\alpha
=\max\limits_{kl\in E}k$ ; b) Minkowski's function $k\left( E,x\right) =$ $%
\max\limits_{k\frac x{\left\| x\right\| }\in E,k>0}k$ is continuous and
positive, $x\in R^n\backslash \left\{ 0\right\} $. Really, the set $E$ is
bounded. Therefore the proposition a) follows from definition \ref{starry}.
Let us prove b). Consider $x_p$, $x\in R^n\backslash \left\{ 0\right\} $, $%
\lim\limits_{p\rightarrow \infty }x_p=x$, $l_p=\frac{x_p}{\left\|
x_p\right\| }$, $l=\frac x{\left\| x\right\| }$. There exists a unique point 
$z_p\in \partial E$ so that $\frac{z_p}{\left\| z_p\right\| }=l_p$. Inasmuch
as $\partial E$ is compact, then we can choose a subsequence $\left\{
z_m\right\} \subset \left\{ z_p\right\} $, $\lim \limits _{m \to \infty}
z_m=u$ and $u\in \partial E$. Hence $\lim\limits_{p\rightarrow \infty
}l_p=l=\frac u{\left\| u\right\| }$. Taking into consideration that $k\left(
E,x_m\right) =\left\| z_m\right\| $, $k\left( E,x\right) =\left\| u\right\| $
and $\left\| z_m\right\| \rightarrow \left\| u\right\| $, $m\rightarrow
\infty $ we have $k\left( E,x_m\right) \rightarrow k\left( E,x\right) $, $%
m\rightarrow \infty $. This proves b).

\section{Necessary and sufficient condition for practical stability}

Let us consider a system of differential equations 
\begin{equation}
\frac{dx}{dt}=f\left( x,t\right)  \label{sys}
\end{equation}
where $x$ is a $n$-dimensional vector, the vector function $f\left(
x,t\right) $ satisfies the conditions of the existence and uniqueness
theorem, $f\left( 0,t\right) \equiv 0$, $t\in \left[ t_0,T\right] $.

We denote $x\left( t\right) =x\left( t,x_0,t_0\right) $ the trajectory of
the system (\ref{sys}) at Cauchy condition $x\left( t_0\right) =x_0$ and
also assume that $\Phi _t\subset R^n$ are compacts, $0\in \Phi _t$, $t\in
\left[ t_0,T\right] $, $G_0\subseteq \Phi _{t_0}$, $0\in G_0$.

\begin{definition}
The unperturbed solution $x\left( t\right) \equiv 0$ of the system (\ref{sys}%
) is said to be $\{G_0,\Phi _t,t_0,T\}$-stable if $x\left( t,x_0,t_0\right)
\in \Phi _t$, $t\in \left[ t_0,T\right] $ as soon as $x_0\in G_0.$
\end{definition}

\begin{theorem}
\label{ndu}For the trivial solution of the system (\ref{sys}) to be $%
\{G_0,\Phi _t,t_0,T\}$-stable it is necessary and sufficient that there
exists a continuous nonincreasing on the system (\ref{sys}) solutions
Lyapunov function $V\left( x,t\right) $ such that 
\begin{equation}
\left\{ x\in R^n:V\left( x,t\right) \leq 1\right\} \subseteq \Phi _t,t\in
\left[ t_0,T\right] ,  \label{restr}
\end{equation}
\begin{equation}
G_0\subseteq \left\{ x\in R^n:V\left( x,t_0\right) \leq 1\right\} .
\label{restr_0}
\end{equation}
\end{theorem}

{\bf Proof.} Necessity. Let us $G_{*}$ be the optimal by inclusion set for $%
\{G_0,\Phi _t,t_0,T\}$-stability of the unperturbed solution $x\left(
t\right) \equiv 0$. It means that zero solution of the system (\ref{sys}) is 
$\{G_{*},\Phi _t,t_0,T\}$-stable and if the solution $x\left( t\right)
\equiv 0$ is $\{D_0,\Phi _t,t_0,T\}$-stable then $D_0\subseteq G_{*}$.
Further, the set $G_{*}$ is compact. Following by lemma \ref{lem} one can
build a function $\alpha \in C\left( R^n\right) $ such that $\alpha (x)=1$
at $x\in \partial G_{*}$, $\alpha \left( x\right) <1$ at $x\in intG_{*}$ and 
$\alpha \left( x\right) >1$ at $x\in R^n\backslash intG_{*}$. Let us
consider the function $V\left( x,t\right) =\alpha \left( \varphi \left(
t,t_0,x\right) \right) $ where $\varphi \left( t,t_0,x\right) =x_0$ on the
solution $x\left( t,x_0,t_0\right) $. Since $V\left( x\left(
t,x_0,t_0\right) ,t\right) =\alpha \left( x_0\right) $ the function $V\left(
x,t\right) $ is nonincreasing on the system (\ref{sys}) solutions. We shall
prove (\ref{restr}) by contradiction. Let there are $x\in R^n$ and $t\in
\left[ t_0,T\right] $ for which the condition (\ref{restr}) is true, but $%
x\notin \Phi _t$. Then $V\left( x,t\right) =V\left( \alpha \left( \varphi
\left( t,t_0,x\right) \right) \right) =\alpha \left( x_0\right) \leq 1$ and $%
x_0\in G_{*}$. Therefore $x=x\left( t,x_0,t_0\right) \in \Phi _t$. We obtain
a contradiction. The validity of (\ref{restr_0}) follows from relations $%
G_0\subseteq G_{*}=\left\{ x\in R^n:\alpha \left( x\right) \leq 1\right\}
=\left\{ x\in R^n:V\left( x,t_0\right) \leq 1\right\} $.

Sufficiency. Suppose that conditions (\ref{restr}) and (\ref{restr_0}) are
satisfied but there are $\tau \in \left[ t_0,T\right] $ and $x_0\in \Phi
_{t_0}$ such that the trajectory $x\left( t\right) =x\left( t,x_0,t_0\right) 
$ leaves the set $\Phi _\tau $. From (\ref{restr}) follows inequality $%
V\left( x\left( \tau \right) ,\tau \right) >1$ which contradicts the
condition (\ref{restr_0}). Our assumption is false and the assertion of the
theorem is true.

\begin{corollary}
\label{V=}If the trivial solution of the system (\ref{sys}) is $\{G_0,\Phi
_t,t_0,T\}$-stable then the function 
\[
V\left( x,t\right) =\left\{ 
\begin{array}{c}
1+\min\limits_{y\in \partial G_{*}}\rho \left( \varphi \left( t,t_0,x\right)
,y\right) ,x\in \left\{ x\in R^n:\varphi \left( t,t_0,x\right) \in
R^n\backslash G_{*}\right\} , \\ 
1-\min\limits_{y\in \partial G_{*}}\rho \left( \varphi \left( t,t_0,x\right)
,y\right) ,x\in \left\{ x\in R^n:\varphi \left( t,t_0,x\right) \in
G_{*}\right\}
\end{array}
\right. 
\]
satisfies all the theorem \ref{ndu} conditions. Here $\rho $ is a metric
equivalent to Euclidean one. If $G_0$ is compact then it is possible to take
the Lyapunov function in the form 
\[
V\left( x,t\right) =\left\{ 
\begin{array}{c}
1+\min\limits_{y\in \partial G_0}\rho \left( \varphi \left( t,t_0,x\right)
,y\right) ,x\in \left\{ x\in R^n:\varphi \left( t,t_0,x\right) \in
R^n\backslash G_0\right\} , \\ 
1-\min\limits_{y\in \partial G_0}\rho \left( \varphi \left( t,t_0,x\right)
,y\right) ,x\in \left\{ x\in R^n:\varphi \left( t,t_0,x\right) \in
G_0\right\} .
\end{array}
\right. 
\]
\end{corollary}

\section{Starry sets of initial conditions: criterion for practical
stability of dynamical system}

We assume the set $G_0\subset R^n$ is a starry compact. Let us introduce a
new variable $y=\left\{ 
\begin{array}{c}
\frac x{k\left( G_0,x_0\right) },x_0\neq 0, \\ 
0,x_0=0
\end{array}
\right. $ and denote $y_0=\left\{ 
\begin{array}{c}
\frac{x_0}{k\left( G_0,x_0\right) },x_0\neq 0, \\ 
0,x_0=0.
\end{array}
\right. $ If $x_0\neq 0$ then 
\[
\frac{dy}{dt}=\frac{dx}{dt}k^{-1}\left( G_0,x_0\right) =f\left( k\left(
G_0,x_0\right) y,t\right) k^{-1}\left( G_0,x_0\right) = 
\]
\[
=f\left( k\left( G_0,y_0\right) y,t\right) k^{-1}\left( G_0,y_0\right) . 
\]

Putting $g\left( y,t,y_0\right) =\left\{ 
\begin{array}{c}
f\left( k\left( G_0,y_0\right) y,t\right) k^{-1}\left( G_0,y_0\right)
,y_0\neq 0, \\ 
0,y_0=0
\end{array}
\right. $ we obtain the system of differential equations 
\begin{equation}
\frac{dy}{dt}=g\left( y,t,y_0\right) ,t\in \left[ t_0,T\right]  \label{g_sys}
\end{equation}
with initial condition $y\left( t_0\right) =y_0\in K_1\left( 0\right) $.

Thus, the problem of $\{G_0,\Phi _t,t_0,T\}$-stability of the system (\ref
{sys}) zero solution can be formulated as the problem of $\{K_1\left(
0\right) ,\Psi _t,t_0,T\}$-stability of the system (\ref{g_sys}) trivial
solution, where $\Psi _t=\{\left( y,y_0\right) :yk\left( G_0,y_0\right) \in
\Phi _t\}$, $y_0\in K_1\left( 0\right) $. From theorem \ref{ndu} follows the
next

\begin{proposition}
For the trivial solution of the system (\ref{g_sys}) to be $\{K_1\left(
0\right) ,\Psi _t,t_0,T\}$-stable it is necessary and sufficient that there
exists a Lyapunov function $V\left( y,t,y_0\right) $ continuous and
increasing on the system (\ref{g_sys}) solutions and 
\[
\left\{ \left( y,y_0\right) :V\left( y,t,y_0\right) \leq 1\right\} \subseteq
\Psi _t,t\in \left[ t_0,T\right] , 
\]
\[
K_1 (0) \subseteq \left\{ y_0:V\left( y_0,t_0,y_0\right) \leq
1\right\} . 
\]
\end{proposition}

If $\rho $ is Euclidean metric then  the Lyapunov function can be taken 
in the form 
\[
V\left( y,t,y_0\right) =\min\limits_{z\in S_1\left( 0\right) }\rho \left(
z,\varphi \left( t,t_0,y,y_0\right) \right) . 
\]
Here $\varphi \left( t,t_0,y,y_0\right) $ is the system (\ref{g_sys}) first
integral in Cauchy's form.

Let us consider the function 
\[
W\left( x,t,x_0\right) =V\left( \frac x{k\left( G_0,x_0\right) },t,\frac{x_0%
}{k\left( G_0,x_0\right) }\right) = 
\]
\[
=\min \limits _{z\in S_1\left( 0\right) } \rho \left( z,\varphi \left(
t,t_0,\frac x{k\left( G_0,x_0\right) },\frac{x_0}{k\left( G_0,x_0\right) }%
\right) \right)=
\]
\[
=\left\| \varphi
\left( t,t_0,\frac x{k\left( G_0,x_0\right) },\frac{x_0}{k\left(
G_0,x_0\right) }\right) \right\| .
\]

\begin{theorem}
The trivial solution of the system (\ref{sys}) is $\{G_0,\Phi _t,t_0,T\}$%
-stable if and only if there exists a continuous Lyapunov function $W\left(
x,t,x_0\right) $ such that $\left( \frac{dW}{dt}\right) _{\left( \ref{sys}%
\right) }\leq 0$ and 
\[
\left\{ \left( x,x_0\right) :W\left( x,t,x_0\right) \leq 1\right\} \subseteq
\left\{ \left( x,x_0\right) :x\in \Phi _t,x_0\in G_0\right\} ,t\in \left[
t_0,T\right] , 
\]
\[
G_0 \subseteq \left\{ x_0:W\left( x_0,t_0,x_0\right) \leq
1\right\} . 
\]
\end{theorem}

\section{Practical stability estimations for linear system}

We assume that the set of initial conditions is a ball, that is $%
G_0=K_c\left( z_0\right) $, and phase restrictions are the sets 
\[
\Gamma _t=\left\{ x\in R^n:\left| l_s^T\left( t\right) \ x\right| \leq
1,s=1,2,\ldots ,N\right\} ,t\in \left[ t_0,T\right] . 
\]

Here $l_s\left( t\right) $ is a continuous $n$-dimensional vector function, $%
s=1,2,\ldots ,N$, parameter $c$ is positive, $z_0\in R^n$.

Following by corollary \ref{V=} Lyapunov function can be chosen in the form 
\[
V\left( x,t\right) =\left\{ 
\begin{array}{c}
1+\min\limits_{y:\left\| y-z_0\right\| =c}\left\| \varphi \left(
t,t_0,x\right) -y\right\| ,x\in \left\{ x:\varphi \left( t,t_0,x\right) \in
R^n\backslash G_0\right\} , \\ 
1-\min\limits_{y:\left\| y-z_0\right\| =c}\left\| \varphi \left(
t,t_0,x\right) -y\right\| ,x\in \left\{ x:\varphi \left( t,t_0,x\right) \in
G_0\right\}
\end{array}
\right. = 
\]
\[
=1+c-\left\| \varphi \left( t,t_0,x\right) -z_0\right\| . 
\]

Now consider a linear system of differential equations 
\begin{equation}
\frac{dx}{dt}=A\left( t\right) x,  \label{l_sys}
\end{equation}
where $A\left( t\right) $ is matrix of order $n\times n$ with continuous
components, $t\in \left[ t_0,T\right] $. In such a case $V\left( x,t\right)
=1-c+\left\| X\left( t_0,t\right) \ x-z_0\right\| $. Here $X\left(
t,t_0\right) $ is a fundamental matrix of solutions of the system (\ref
{l_sys}) normalized with respect to $t_0$. The problem is to define the
maximal radius $c^{*}$ of the boll $K_c\left( z_0\right) $ such that the
trivial solution of the system (\ref{l_sys}) is $\left\{ K_{c^{*}}\left(
z_0\right) ,\Phi _t,t_0,T\right\} $-stable. From theorem \ref{ndu} follows
that we choose parameter $c$ so that $\left| l_s^T\left( t\right) \ x\right|
\leq 1$, $s=1,2,\ldots ,N$ as soon as $V\left( x,t\right) \leq 1$.
Maximizing Lagrange function 
\[
L\left( x,\lambda \right) =l_s^Tx+\lambda \left( \left\| X\left(
t_0,t\right) \ x-z_0\right\| ^2-c^2\right) ,s\in \left\{ 1,2,\ldots
N\right\} , 
\]
$t\in \left[ t_0,T\right] $ yield 
\begin{equation}
c^{*}=\min\limits_{t\in \left[ t_0,T\right] }\min\limits_{s=1,2,\ldots N}%
\frac{1-\left| l_s^T\left( t\right) X\left( t,t_0\right) \ z_0\right| }{%
\sqrt{l_s^T\left( t\right) H\left( t\right) l_s\left( t\right) }},
\label{radius}
\end{equation}
where the matrix $H\left( t\right) =X\left( t,t_0\right) X^T\left(
t,t_0\right) $ satisfies the Lyapunov's differential matrix equation 
\begin{equation}
\frac{dH\left( t\right) }{dt}=A\left( t\right) H\left( t\right) +H\left(
t\right) A^T\left( t\right) ,t\in \left[ t_0,T\right]  \label{L_sys}
\end{equation}
with initial condition $H\left( t_0\right) =I$, $I$ is unit matrix of order $%
n$.

Let $Q$ be a positive definite symmetric matrix, 
\[
\hat{K}_c\left( z_0\right) =\left\{ x\in R^n:\left( x-z_0\right) ^TQ\left(
x-z_0\right) \leq c\right\} . 
\]
Consider the set of initial conditions in the form $G_0=\hat{K}_c\left(
z_0\right) $. It is necessary to obtain the maximal $c^{*}$ from such $c>0$
that the trivial solution of the system (\ref{l_sys}) is $\{\hat{K}_c\left(
z_0\right) ,\Gamma _t,t_0,T\}$-stable. If we choose $\rho \left( a,b\right) =%
\sqrt{\left( a-b\right) ^TQ\left( a-b\right) }=\left\| a-b\right\| _Q$ then $%
V\left( x,t\right) =1-c+\left\| X\left( t_0,t\right) \ x-z_0\right\| _Q$.
Proceeding as it stated above we obtain 
\[
c^{*}=\min\limits_{t\in \left[ t,T\right] }\min\limits_{s=1,2,\ldots N}\frac{%
1-\left| l_s^T\left( t\right) X\left( t,t_0\right) z_0\right| }{\sqrt{%
l_s^T\left( t\right) H\left( t\right) l_s\left( t\right) }}. 
\]
Here $H\left( t\right) =X\left( t,t_0\right) Q^{-1}X^T\left( t,t_0\right) $
satisfies (\ref{L_sys}), $H\left( t_0\right) =Q^{-1}$.

\end{document}